\documentclass[11pt]{amsart}
\usepackage{geometry}
\usepackage[utf8]{inputenc}
\usepackage{graphicx}
\usepackage{amssymb}
\usepackage{epstopdf}
\usepackage{setspace}

\title[Minkowski's Question-Mark functions]{Generalised L\"uroth expansions and a family of Minkowski's Question-Mark functions}

\author{Aubin Arroyo}
\date{\today} 
\address{Instituto de Matem\'aticas, Unidad Cuernavaca, Universidad Nacional Aut\'onoma
de M\'exico, A.P. 273-3 Admon. 3, Cuernavaca, Morelos, 62251, M\'exico.}
\email{aubinarroyo@im.unam.mx} 
\thanks{Partially supported by CONACyT-167594, México.}
\newcommand{\interval}{\mathcal{I}}
\newcommand{\Q}{?}
\newcommand{\cf}{\mbox{cf}\,}
\newcommand{\lur}{\mbox{L\"ur}_\alpha}

\newtheorem{theorem}{Theorem}
\newtheorem{lemma}{Lemma}
\newtheorem{e-definition}{Definition}

\begin{document}

\begin{abstract}
The Minkowski's Question-Mark function is a singular homeomorphism of the unit interval that maps the set of quadratic surds into the rationals.
This function has deserved the attention of several authors since the beginning of the twentieth century. 
Using different representations of real numbers by infinite sequences of integers, called $\alpha$-L\"uroth expansions, we obtain different instances of the standard shift map on infinite symbols, all of them topologically conjugated to the Gauss Map. In this note we prove that each of these conjugations share properties with the Minkowski's Question-Mark function: all of them are singular homeomorphisms of the interval, and in the ``rational'' cases, they map the set of quadratic surds into the set of rational numbers. 
In this sense, this family is a natural generalisation of the Minkowski's Question-Mark function.
\end{abstract}

\maketitle

A \emph{quadratic surd} is a real solution of a quadratic polynomial with integer coefficients. 
Denote by $\mathbb{A}_2$ the set of quadratic surds in the unit interval $\interval := [0,1] \subset \mathbb{R}$, and denote by $\mathbb{Q}$ the set of rational numbers contained in $\interval$. 
Minkowski, in \cite{Minkowski11}, defined a function which is a one to one correspondence between the set $\mathbb{A}_2 \cup \mathbb{Q}$ and the set $\mathbb{Q}$. This function is called \emph{Minkowski's Question-Mark function}, and is denoted by $\Q: \interval \to \interval$. 

On the other hand, the function $\Q$  is a singular function, that is,
its derivative exists and~$\Q'(x) = 0$ for Lebesgue almost every $x \in \interval$. 
Functions like Cantor maps, sometimes called \emph{devil staircases}, are singular maps if they are defined using a Lebesgue measure zero Cantor set on the interval; however these functions are constant on each connected component of the complement of such Cantor set. The function $\Q$ is strictly increasing and singular, so it is often called a \emph{slippery devil staircase}.

The Minkowski's Question-Mark function is  defined using Farey fractions and mediants: set $\Q(0/1) = 0$, $\Q(1/1)=1$ and define it in the following way:  if $p/q$ and $p'/q'$ are consecutive Farey rationals of order $n>1$, set:
$$\Q\left(\frac{ p+p'}{q+q'}\right) = \frac{\Q(p/q) + \Q(p'/q')}{2}.$$ 
After Minkowski several generalisations have been constructed, beginning  with Denjoy and Salem; see  \cite{Denjoy32}, \cite{Denjoy38} and \cite{Salem43}. Much more recently, Kessebohmer, Munday and Stratmann in \cite{KessebohmerMundayStratmann12}, and Munday  in \cite{Munday14}, study this map from a dynamical systems point of view, and the relation with the generalised L\"uroth expansion of numbers. In this note we exploit further this relation to obtain a family of functions which are generalisations of Minkowski's Question Mark function. 

Let us recall the continued fraction expression of a number. Given $x \in\interval$, the continued fraction expansion of $x$ is a sequence of positive integers $\cf(x):=[a_1,a_2, \ldots] \in \mathbb{N}^{\mathbb{N}}$ such that:
$$x = \frac{1}{a_1 + \frac{1}{a_2 +\cdots }}.$$
For rational numbers this sequence is finite and is not unique. In fact, $[a_1,\ldots, a_n] = [a_1,\ldots, a_n-1,1]$. 
Conversely, if $x \in \interval \smallsetminus \mathbb{Q}$ then $\cf(x)$ is infinite and unique.
Given $x \in \interval$ such that $\cf(x)=[a_1,a_2, \ldots]$, the formula given in \cite{Salem43} is:  
\begin{equation} \label{Salem}
\Q (x) = \frac{1}{2^{a_1-1}} -\frac{1}{2^{a_1+a_2-1}} + \frac{1}{2^{a_1+a_2+a_3-1}} - \cdots 
\end{equation}
Notice that for rational points $x\in \mathbb{Q}$, the function $\Q(x)$ is well defined and its value is rational, since the formula involves only a finite number of terms.

Generalised $\alpha$-L\"uroth expansions are  representations of real numbers by a sequence of integers, depending on a partition of the interval.
They where studied in \cite{Barrionuevo96} and are a generalisation of the classical \emph{L\"uroth expansion} of \cite{Luroth83}.
In this note, we will restrict our attention to the class of $\alpha$-L\"uroth expansions which were introduced in  \cite{KessebohmerMundayStratmann12} (see also  \cite{Munday14}).
Each of these representations depends on a countable partition of $\interval$ defined by a sequence with the following properties:
Let $\alpha = \{t_j\}_{j\geq 1} \subset \interval$
be a strictly decreasing sequence such that $t_1 = 1$, and $t_j \to 0$ when $j\to \infty$. 
This sequence defines a partition of $\interval$ by non-trivial intervals:
\begin{equation}\label{parti} 
\left\{ (t_{j+1}, t_{j}] \,|\, j \geq 1\right\}
\end{equation}
of length  $t_j-t_{j+1}>0$. 
In this way, for any $x \in \interval$, there is a sequence of positive integers $\lur(x) := (a_n)_{n\geq 1}$ (finite or infinite) such that: 
\begin{equation} \label{family}
x = t_{a_1} + \sum_{n=2}^{\infty}  (-1)^{n+1}  \prod_{j=1}^{n-1}(t_{a_j}-t_{a_j+1})\, t_{a_n}. 
\end{equation}
In this way, $\lur(x)$ is called 
the $\alpha$-L\"uroth expansion of $x$ with respect to the given partition $\alpha$.
This expression can be found in \cite{KessebohmerMundayStratmann12} and \cite{Munday14}. 
Moreover, except for a countable set in $\interval$, this expression is unique.
It is not difficult to see that the formula in \eqref{Salem} can be derived from \eqref{family} setting $t_n = 2^{n-1}$, for $n\geq1$. On the other hand, if we set $t_n = 1/n$, we obtain in the standard alternating L\"uroth expansion of $x$.

We will show that for any sequence $\alpha$ there is a singular homeomorphism of the interval~$\Q_\alpha: \interval \to \interval$. Moreover, if $\alpha \subset \mathbb{Q}$ then $\Q_\alpha$ maps quadratic surds into the set of rational numbers. 
In words, to define $\Q_\alpha(x)$ for $x\in \interval$ one needs to compute the continued fraction expression of $x$. The real number which has this sequence of integers as its generalised $\alpha$-L\"uroth expansion is $\Q_\alpha(x)$. 
The purpose of this note is to prove the following theorem, which states that these maps are  generalisations of the Minkowski's Question-Mark function.

\begin{theorem} \label{elmero}
Let $\alpha$ be a sequence as in (\ref{parti}). The function:
$$\Q_{\alpha} (x) := t_{a_1} + \sum_{n=2}^{\infty}  (-1)^{n+1}    \prod_{j=1}^{n-1}(t_{a_j}-t_{a_j+1}) \, t_{a_n},$$
where $\cf(x) = [a_1, a_2, \ldots]$,
is a  singular homeomorphism of $\interval$. Moreover, if $\alpha \subset \mathbb{Q}$, then~$\Q_{\alpha}( \mathbb{A}_2 ) \subset \mathbb{Q}$.  
\end{theorem}

The proof of the fact that these functions are singular is a consequence of the fact that ergodic invariant measures for dynamical system either coincide or are mutually singular.

\section{Proof of results}
Let $\Sigma := \mathbb{N}^\mathbb{N}$ be the set of infinite sequences of positive integers with the product topology. Let $\sigma: \Sigma \to \Sigma$ be the shift map, that is $\sigma(\, (a_1,a_2,\ldots) \,) = (a_2,a_3,\ldots)$. 
On the other hand, consider the Gauss map $G:\interval \to \interval$ defined by: 
$$G(x)=\frac{1}{x}  -   \left \lfloor \frac{1}{x} \right \rfloor, $$
for $x \neq 0$, and $G(0)=0$.
It is a well known fact that the set $\mathbb{A}_2$ corresponds to the set of pre-periodic points of $G$, $\mbox{Per}(G) := \{y \in \interval \,| \,\exists n,m \in \mathbb{N} \mbox{ such that } G^{n+m}(x)  = G^n(x)\}$.
Furthermore, the Gauss map is a dynamical factor of $\sigma$. 
In fact, the inverse map $\cf^{-1}:  \Sigma \to \interval$
is a continuous and surjective map that satisfies 
$\sigma \circ \cf^{-1}=  \cf^{-1} \circ G $. 
If we set the partition $\gamma$ of $\interval$, given by $I_n = ( (n+1)^{-1}, n^{-1}]$, for $n\geq 1$, we have that if $\cf(x) = (a_1, a_2, \ldots)$, then $G^k(x) \in I_{a_k}$, for any $k\geq 1$.
That is, the continued fraction expression of a point $x$ describes the itinerary $x$ under the iterations of $G$, with respect to the partition $\gamma$.

The following class of maps, introduced in  \cite{KessebohmerMundayStratmann12},  are also factors of the shift map~$\sigma$. 
\begin{e-definition}
Given a sequence $\alpha$ as in (\ref{parti}), an $\alpha$-L\"uroth map is a piecewise linear map
$L_{\alpha}: \interval \to \interval$ defined by $L_{\alpha}(0)=0$ and:
\begin{equation}
L_{\alpha}(x) = \left \{ \frac{ t_j\,-\,x}{  t_j-t_{j+1}} , \mbox{ if } x \in  [t_{j+1}, t_{j}], \, j \in \mathbb{N}. \right .
\end{equation}
\end{e-definition}

For these maps, the $\alpha$-L\"uroth expansion of  $x$ codifies the itinerary of the orbit of $x$, under the iterations of $L_\alpha$, with respect to the partition $\alpha$. In fact,  we have that $\lur \circ \sigma = L_\alpha \circ \lur $. 
Precisely, if $\lur(x) := (a_1, a_2, \ldots)$ then $L_\alpha^k(x) \in (t_{a_k+1}, t_{a_k}]$, for any $k\geq 1$.
Moreover, if $\alpha \subset \mathbb{Q}$, then the set of pre-periodic points is $\mathbb{Q}$.

\subsection{Invariant measures}Denote by $\lambda$ the Lebesgue measure on $\interval$.
The Gauss measure, given by $\mu_G(A) := \frac{1}{\log 2}\int_A\frac{dx}{1+x}$, for any Borel set $A\subset \interval$, is the unique invariant measure for $G$ which is absolutely continuous with respect to $\lambda$. 
On the other hand, the Lebesgue measure is the unique ergodic invariant measure which is absolutely continuous with respect to $\lambda$, for any $\alpha$-L\"uroth map, for any given sequence $\alpha$.
The following lemma, whose proof we omit, is a consequence of a stronger statement that can be found in \cite{KessebohmerMundayStratmann12}.
The uniqueness of this measure can also be obtained from the Folklore Theorem of existence and uniqueness of absolutely continuous invariant measures for Markov Maps; see \cite{Bowen79}.
 
\begin{lemma} \label{two:lemB}
If $T$ is a generalised linear $\alpha$-L\"uroth map, then $\lambda$, the Lebesgue measure, is the unique $T$-invariant and ergodic measure which is absolutely continuous with respect to the Lebesgue measure. 
\end{lemma}

Now we are in a position to state and prove the following theorem, which proves Theorem~\ref{elmero}. 

\begin{theorem} \label{teo:main}
Given a partition $\alpha$ as above, the map $\Q_\alpha:\interval \to \interval$ is a increasing homeomorphism such that $\Q_\alpha$ is singular and that $L_\alpha \circ \Q_\alpha = \Q_\alpha \circ G$.
Moreover,  if $\alpha \subset \mathbb{Q}$ then $\Q_\alpha(\mathbb{A}_2) \subset \mathbb{Q}$.
\end{theorem}

\noindent \emph{Proof.}
Let $\alpha$ be a partition as above. 
If we define
$\Q_\alpha := \lur ^{-1}\circ \cf : \interval \to \interval$, then clearly $ L_\alpha \circ \Q_{\alpha} = \Q_{\alpha} \circ G$, and hence is a topological conjugation between $G$ and $L_\alpha$. 
Notice that~$\Q_\alpha'(x)$ exists for $\lambda$-a.e $x \in \interval$, since the function $\Q_\alpha$ is increasing.
Moreover, if $\alpha \subset \mathbb{Q}$ then $\Q_\alpha(\mathbb{A}_2) = \mbox{Per}(L_\alpha) \subset \mathbb{Q}$, since $\Q_\alpha$ is a conjugation.

To prove that $\Q_\alpha$ is singular we need to show that there is a Borel set $\tilde B$ such that $\lambda(\tilde B)=0$ and $\lambda  (\Q_\alpha(\tilde B)) =1$; see Theorem 3.72 of \cite{Giovanni09}, for instance.
For that, consider the pull-back measure $\mu_\alpha(B) := \lambda(\Q_\alpha(B))$, for any Borel set $B$ on $\interval$. 
The measure $\mu_\alpha$ is invariant and ergodic for $G$.
This follows from the fact that $\lambda$ is an invariant ergodic measure for $L_\alpha$. 

So, we have two ergodic measures for $G$. Therefore, 
either $\mu_\alpha = \mu_G$ or they are singular to each other. If $\mu_\alpha = \mu_G$, we have that for any Borel set $B$: 
$$ \int_{\Q_\alpha(B)} dx =   \int_B  \Q_\alpha^{-1}(x) dx = \frac{1}{\log \, 2} \int_B \frac{dx}{1+x}. $$
Therefore,  $\Q_\alpha^{-1}(x) = 1/( \log \, 2(1+x))$, for $\lambda$-a.e $x\in \interval$. Observe that the right side hand is a decreasing function on $x$, and that $\Q_\alpha^{-1}(x)$ is increasing. This is a contradiction. Then,  $\mu_\alpha$ and $ \mu_G$ are mutually singular. Hence, there exists a measurable set $\tilde B \subset \interval$ such that $\lambda(\tilde B)=0$ and $\lambda  (\Q_\alpha(\tilde B)) =1$, and we are done.
\qed

\section*{Acknowledgements}
The author  is grateful to Carlos Cabrera (UCIM-UNAM)
and Sara Munday (University of York) for many helpful discussions on this subject.  
This work was partially supported by CONACyT-167594, México.

\end{document}